\documentclass[12pt]{amsart}
\usepackage{multicol}
\usepackage[numbers,sort&compress]{natbib}
\usepackage[euler-digits]{eulervm}
\usepackage[a4paper,body={14.9cm,23cm},centering,includeheadfoot]{geometry} 
\usepackage{amssymb}
\usepackage{xy}
\usepackage{ctable}
\usepackage{longtable}
\usepackage{array}
\usepackage{aliascnt}
\usepackage{breakurl}
\usepackage[noend]{algpseudocode}
\usepackage{algorithm}
\usepackage[colorlinks, breaklinks, linkcolor=blue]{hyperref} 

\newtheorem{theorem}{Theorem}
\def\equationautorefname~#1\null{(#1)\null}
\def\itemautorefname~#1\null{(#1)\null}
\def\sectionautorefname~#1\null{\S#1\null}
\def\subsectionautorefname~#1\null{\S#1\null}

\newcommand{\GAP}{\textsf{GAP}} 
\newcommand{\CHEVIE}{\textsf{CHEVIE}}
\newcommand{\ZigZag}{\textsf{ZigZag}} 
\numberwithin{equation}{section}

\newcommand{\otoprule}{\midrule[\heavyrulewidth]}
\newcommand\Ind{\operatorname{Ind}}

\title[Computations for rank five and six Coxeter groups]
{Computations for Coxeter arrangements and Solomon's descent algebra II:
  Groups of rank five and six}

\keywords{Coxeter group, descent algebra, Orlik-Solomon algebra}
\subjclass[2010]{20F55, 20C15, 20C40, 52C35}

\begin{document}

\allowdisplaybreaks 
\author[M. Bishop]{Marcus Bishop} 

\address{Fakult\"at f\"ur Mathematik\\Ruhr-Universit\"at Bochum\\D-44780
  Bochum, Germany} 

\email{marcus.bishop@rub.de}

\author[J. M. Douglass]{J. Matthew Douglass} 

\address{Department of Mathematics\\University of North Texas\\Denton TX,
  USA 76203} 

\email{douglass@unt.edu} 

\author[G. Pfeiffer]{G\"otz Pfeiffer} 

\address{School of Mathematics, Statistics and Applied Mathematics\\National
  University of Ireland, Galway\\University Road, Galway, Ireland}

\email{goetz.pfeiffer@nuigalway.ie} 

\author[G. R\"ohrle]{Gerhard R\"ohrle}

\address{Fakult\"at f\"ur Mathematik\\Ruhr-Universit\"at Bochum\\D-44780
  Bochum, Germany} 

\email{gerhard.roehrle@rub.de}

\begin{abstract}
  In recent papers we have refined a conjecture of Lehrer and Solomon
  expressing the character of a finite Coxeter group $W$ acting on the
  graded components of its Orlik-Solomon algebra as a sum of
  characters induced from linear characters of centralizers of elements of
  $W$. The refined conjecture relates the character above to a decomposition
  of the regular character of $W$ related to Solomon's descent algebra of
  $W$. The refined conjecture has been proved for symmetric and dihedral
  groups, as well as for finite Coxeter groups of rank three and four.  In
  this paper, we prove the conjecture for finite Coxeter groups of rank
  five and six. The techniques developed and implemented in this
  paper provide previously unknown decompositions of the regular and
  Orlik-Solomon characters of the groups considered.
\end{abstract}

\maketitle

\section{Introduction}
Let $\left(W,S\right)$ be a finite Coxeter system.  In previous articles
\cite{rank3and4,douglasspfeifferroehrle:inductive,
  douglasspfeifferroehrle:coxeter} we proposed a conjecture relating the
character $\omega$ of the Orlik-Solomon algebra of $W$ to the regular
character $\rho$ of $W$. 
Based on a conjecture of Lehrer and Solomon~\cite{lehrersolomon:symmetric}, 
in this paper we prove the conjecture for the
Coxeter groups of type $B_5$, $B_6$, $D_5$, $D_6$, and $E_6$.
These computations, together with the remarks about reducible Coxeter
groups following Theorem~2.3 of \cite{rank3and4} and the proof of the
conjecture for groups of type $A$ in \cite{douglasspfeifferroehrle:coxeter},
prove the conjecture for all finite Coxeter groups of rank five and six.
Our result is stated for these groups as the following theorem.

\begin{theorem}\label{ConjectureA}
  Suppose that $W$ is a finite Coxeter group of rank five or six and
  that $\mathcal{R}$ is a set of conjugacy class representatives of
  $W$. Then for each $w\in\mathcal{R}$ there exists a linear character
  $\varphi_w$ of $C_W(w)$ such that
  \[
  \rho= \sum_{w\in\mathcal{R}} \mathrm{Ind}_{C_W(w)}^W \varphi_w
  \quad\text{and}\quad\omega= \epsilon \sum_{w\in\mathcal{R}}
  \mathrm{Ind}_{C_W(w)}^W (\alpha_w \varphi_w )
  \] 
  where $\epsilon$ is the sign character of $W$ and for $z\in C_W(w)$,
  $\alpha_w(z)$ denotes the determinant of the restriction of $z$ to the 
  $1$-eigenspace of $w$ in the complex reflection representation of $W$.
\end{theorem}

Let $A\left(W\right)$ be the Orlik-Solomon algebra of $W$.  The strategy for
proving~\autoref{ConjectureA} is to decompose the $\mathbb{C}W$-modules
$\mathbb{C}W$ and $A\left(W\right)$ into direct sums and prove a refinement
of~\autoref{ConjectureA} for each summand.  This method is somewhat stronger
than directly proving~\autoref{ConjectureA}, because it requires the
solution to be compatible with the direct sum decompositions of
$\mathbb{C}W$ and $A\left(W\right)$. This method also has the advantage that
it splits the problem into smaller problems and provides additional insight
into how the representations of $W$ on $\mathbb{C}W$ and on
$A\left(W\right)$ are related.

The decomposition of $\mathbb{C}W$ comes from idempotents $e_\lambda$ in the
descent algebra $\Sigma\left(W\right)$ constructed in
\cite{bergeronbergeronhowletttaylor:decomposition}.  These idempotents are
indexed by subsets of $S$ up to conjugacy in $W$. A class of conjugate
subsets of $S$ is called a {\em shape} of $W$. Denote the set of shapes of
$W$ by $\Lambda$. In \cite{bergeronbergeronhowletttaylor:decomposition} it
is shown how to construct a quasi-idempotent $e_L$ for any $L\subseteq S$
and then $e_\lambda$ is the sum of the quasi-idempotents $e_L$ where $L$
runs over the subsets in the shape $\lambda$.  It is also shown that
$\left\{ e_\lambda \mid \lambda\in\Lambda \right\}$ is a complete set of
primitive orthogonal idempotents of $\Sigma\left(W\right)$.  Since
$\Sigma\left(W\right)$ is a subalgebra of $\mathbb{C}W$ and
$1=\sum_{\lambda\in\Lambda}e_\lambda$, we conclude that $\mathbb{C}W=
\bigoplus_{\lambda\in\Lambda} e_\lambda\mathbb{C}W$ as a
$\mathbb{C}W$-module.  Denoting the character of $e_\lambda\mathbb{C}W$ by
$\rho_\lambda$ we have
\begin{equation}
  \label{eq:1a}
  \rho=\sum_{\lambda\in\Lambda}\rho_\lambda.  
\end{equation}

The corresponding decomposition of $A\left(W\right)$ comes from Brieskorn's
Lemma. Let $T$ be the set of reflections in the complex reflection
representation $V$ of $W$. Recall that the Orlik-Solomon algebra
$A\left(W\right)$ may be defined as the quotient of the exterior algebra
with generators $\left\{e_t\mid t\in T\right\}$ by the ideal generated by
elements of the form $\sum_{i=1}^k \left(-1\right)^ie_{t_1} e_{t_2}\cdots
\widehat{e_{t_i}} \cdots e_{t_k}$ for all sets $\left\{t_1,t_2,\ldots,t_k
\right\}\subseteq T$ of linearly dependent reflections.  Here, we say that a
set of reflections is linearly dependent if the linear forms defining their
reflecting hyperplanes are linearly dependent in the dual of $V$.  For $t\in
T$ we denote the image of the generator $e_t$ in $A\left(W\right)$ by
$a_t$. Thus, an arbitrary element of $A\left(W\right)$ can be expressed as a
linear combination of monomials $a_{t_1}a_{t_2}\cdots a_{t_k}$ with
$t_1,t_2,\ldots,t_k\in T$. The algebra $A\left(W\right)$ is a right
$\mathbb{C}W$-module with $\left(a_{t_1}a_{t_2}\cdots a_{t_k}\right).w
=a_{w^{-1}t_1w}a_{w^{-1}t_2w}\cdots a_{w^{-1}t_kw}$ for $w\in W$.

Each monomial $a_{t_1}a_{t_2}\cdots a_{t_k}$ determines a subspace of $V$,
namely the intersection of the fixed point spaces of the reflections
$t_1,t_2,\ldots,t_k$.  If $X$ is a subspace of $V$, then we denote by $A_X$
the span of all monomials with fixed point space equal to $X$.  Then taking
$A_\lambda$ to be the sum of the $A_X$ for which $X$ is the fixed point set
of some conjugate of $W_L$ for some $L\in\lambda$, we have a decomposition
$A\left(W\right) =\bigoplus_{\lambda\in\Lambda}A_\lambda$.  Denoting the
character of $A_\lambda$ by $\omega_\lambda$ we have
\begin{equation*}
  \omega=\sum_{\lambda\in\Lambda}\omega_\lambda.
\end{equation*}

Finally, we choose a set of conjugacy class representatives compatible with
the decompositions above. A conjugacy class in $W_L$ is called {\em
  cuspidal} if the fixed point set in the reflection representation of $W_L$
of any of its elements is trivial.  Now if we choose a fixed representative
$L\left(\lambda\right)$ of each shape $\lambda$ and let
$\mathcal{C}_{L\left(\lambda\right)}$ be a set of representatives of the
cuspidal classes in $W_{L\left(\lambda\right)}$, then by Theorem~3.2.12
of~\cite{geckpfeiffer:characters}
\begin{equation}
  \label{eq:2a}
  \mathcal{R}=\bigcup_{\lambda\in\Lambda}
  \mathcal{C}_{L\left(\lambda\right)} \quad\text{is a set of conjugacy class
    representatives of $W$.}    
\end{equation}

Suppose that $L\subseteq S$. The homogeneous component of $A\left( W_L
\right)$ of highest degree is called the {\em top component} of
$A\left(W_L\right)$.  On the other hand, $W_L$ is also a Coxeter group and
thereby admits a system of quasi-idempotents as in
\cite{bergeronbergeronhowletttaylor:decomposition}, now denoted by $e^L_J$
for $J\subseteq L$ to distinguish them from the quasi-idempotents $e_J$ in
$\mathbb{C}W$. Note that this notation does not agree with that used in
\cite[\S7]{bergeronbergeronhowletttaylor:decomposition}. In analogy with
$A\left(W_L\right)$ the homogeneous component $e_L^L\mathbb{C}W_L$ of
$\mathbb{C}W_L$ is called the {\em top component} of $\mathbb{C}W_L$. We
denote the characters of $W_L$ afforded by the top components of
$A\left(W_L\right)$ and $\mathbb{C}W_L$ by $\omega_L$ and $\rho_L$
respectively.

Now the top components of $A\left(W_L\right)$ and $\mathbb{C}W_L$ are
naturally $W_L$-stable subspaces of $A\left(W\right)$ and $\mathbb{C}W$ and
it turns out that they are $N_W\left(W_L\right)$-stable by
\cite[Proposition~4.8]{douglasspfeifferroehrle:coxeter}.  Thus they afford
characters $\widetilde{\omega_L}$ and $\widetilde{\rho_L}$ of
$N_W\left(W_L\right)$ which are extensions of $\omega_L$ and
$\rho_L$. Furthermore, if $L\in \lambda$ then by the same proposition,
\begin{equation}
  \label{eq:3a}
  \omega_\lambda=\mathrm{Ind}_{N_W\left(W_L\right)}^W
  \widetilde{\omega_L}\quad\text{and} \quad\rho_\lambda=
  \mathrm{Ind}_{N_W\left(W_L\right)}^W \widetilde{\rho_L}.
\end{equation}

Consider the following refinement of \autoref{ConjectureA}. In the statement
of this theorem we use the fact, proved in \cite[Theorem
3.1]{konvalinkapfeifferroever:centralizers}, that if $w$ is cuspidal in
$W_L$, then $C_W\left(w\right)\subseteq N_W\left(W_L\right)$.

\begin{theorem}\label{ConjectureC}
  Suppose that $\left(W,S\right)$ is a finite Coxeter system of rank five or
  six, that $L\subseteq S$, and that $\mathcal{C}_L$ is a set of
  representatives of the cuspidal conjugacy classes of $W_L$. Then for each
  $w\in\mathcal{C}_L$ there exists a linear character $\varphi_w$ of
  $C_W\left(w\right)$ such that
  \[
  \widetilde{\rho_L}=\sum_{w\in\mathcal{C}_L}
  \mathrm{Ind}_{C_W(w)}^{N_W(W_L)} \varphi_w =\alpha_L\epsilon
  \widetilde{\omega_L}
  \]
  where for $n\in N_W(W_L)$, $\alpha_L(n)$ denotes the determinant of
  the restriction of $n$ to the subspace of fixed points of $W_L$ in
  $V$.
\end{theorem}

To prove \autoref{ConjectureA} we prove \autoref{ConjectureC} for the
representative $L\left(\lambda\right)$ of each shape $\lambda$. Then the
characters $\varphi_w$ that satisfy \autoref{ConjectureC} with
$L=L\left(\lambda\right)$ as $\lambda$ varies over all shapes prove
the first equality of \autoref{ConjectureA} because
\begin{align*}
  \rho&=\sum_{\lambda\in\Lambda}\rho_\lambda&\text{by (\ref{eq:1a})}\\
  &=\sum_{\lambda\in\Lambda}\Ind_{N_W\left(W_{L\left(\lambda\right)}\right)}^W
  \widetilde{\rho_{L\left(\lambda\right)}} &\text{by (\ref{eq:3a})}\\
  &=\sum_{\lambda\in\Lambda}\sum_{w\in\mathcal{C}_{L\left(\lambda\right)}}
  \Ind_{N_W\left(W_{L\left(\lambda\right)}\right)}^W
  \Ind_{C_W\left(w\right)}^{N_W\left(W_{L\left(\lambda\right)}\right)}
  \varphi_w&\text{by~\autoref{ConjectureC}}\\
  &=\sum_{w\in\mathcal{R}}\Ind_{C_W\left(w\right)}^W\varphi_w &
\end{align*}
where the last equality follows from transitivity of induction and
(\ref{eq:2a}). A similar argument proves the second equality in
\autoref{ConjectureA}.  We prove \autoref{ConjectureC} in
\autoref{ConjectureBSection} and~\autoref{ConjectureCSection}.

\section{Implementation}
As in \cite{rank3and4}, we have implemented the calculations for this
article in the computer algebra system \GAP\ \cite{gap3} in conjunction with
the \CHEVIE\ \cite{chevie} and the \ZigZag\ \cite{zigzag} packages.  In
addition to our comments about the implementation in \cite{rank3and4} we
make the following remarks about the techniques new to this paper and
improvements to old techniques.

\subsection{The Extension \texorpdfstring{$\widetilde{\rho_L}$}{rhoLtilde}}
\label{RhoWiggleSection}
In this subsection we develop a formula for the character
$\widetilde{\rho_L}$ of $N_W\left(W_L\right)$ for $L\subseteq S$.  First we
review the definitions of the constructions used in the process.

If $J\subseteq L$ then the {\em parabolic transversal} of $W_J$ in $W_L$ is
the set $X_J^L$ of elements $w\in W_L$ satisfying
$\ell\left(sw\right)>\ell\left(w\right)$ for all $s\in J$, where $\ell$ is
the usual length function of $W$ with respect to $S$.  Then $X_J^L$ can be
calculated directly from the definition or by using the {\ttfamily
  ParabolicTransversal} function supplied by the \ZigZag\ package.

In order to use the formula for $\widetilde{\rho_L}$ below, we need to be
able to decompose an element of $N_W\left(W_L\right)$ into the product of an
element of $W_L$ and an element of the normalizer complement $N_L$ of $W_L$.
Recall that $N_L$ consists of certain elements of the parabolic transversal
$X^S_L$ of $W_L$ in $W$.  Therefore, the decomposition of an element of
$N_W\left(W_L\right)$ into a product $nw$ with $n\in N_L$ and $w\in W_L$ is
a special case of the more general decomposition of an element of $W$ into a
product of a coset representative in $X^S_L$ by an element of $W_L$.  In
\ZigZag\ this decomposition is implemented as the {\ttfamily
  ParabolicCoordinates} function.

The quasi-idempotents $e_J^L$ are defined in \cite
{bergeronbergeronhowletttaylor:decomposition} by means of the matrix
$M=\left(m_{KJ}\right)$, whose rows and columns are indexed by the subsets of
$S$ and whose $\left(K,J\right)$-entry is
\[
m_{KJ}=\begin{cases} \left|\left\{x\in X_K\mid J^x\subseteq S\right\}\right|
  &\text{if $K\supseteq J$}\\0&\text{otherwise.}\end{cases}
\] 
The matrix $M$ can be calculated directly from the definition or by calling
the method {\ttfamily Mu} supplied by the \ZigZag\ package.  Then putting
$N=M^{-1}= \left(n_{KJ}\right)$ we have $e_J^L=\sum_{K\subseteq
  L}n_{JK}x_K^L$ where $x_K^L$ is the sum in $\mathbb{C}W_L$ of the elements
of $X_K^L$.

In the following discussion let $w\in W_L$, $n\in N_L$, and
$x\in\mathbb{C}W_L$.  Observe that $N_W\left(W_L\right)$ acts on
$\mathbb{C}W_L$ on the right by $x.\left(wn\right)=n^{-1}xwn$.  Using this
action we define the map
\[
\gamma\left(wn,x\right):\mathbb{C}W_L\to\mathbb{C}W_L\qquad\text{by}\qquad
\gamma\left(wn,x\right)\left(v\right)
=\left(xv\right).\left(wn\right)=n^{-1}xvwn
\] 
for $v\in\mathbb{C}W_L$.

The idempotent $e_L^L\in\mathbb{C}W_L$ determines a
$\mathbb{C}N_W\left(W_L\right)$-stable decomposition 
of the group algebra of $W_L$, $\mathbb{C}W_L
=e_L^L\mathbb{C}W_L\oplus \left(1-e_L^L\right) \mathbb{C}W_L$.  Calculating
the trace of the action of $\gamma\left(wn,e_L^L\right)$ with respect to a basis of
$\mathbb{C}W_L$ adapted to this decomposition, we find that $\widetilde{\rho_L}
\left( wn \right) =\mathsf{Tr} \left( \gamma \left( wn,e_L^L \right)
\right)$ since $\gamma\left(wn,e_L^L\right)$ sends
$\left(1-e_L^L\right)\mathbb{C}W_L$ to $e_L^L\mathbb{C}W_L$.
Using the linearity of $\gamma$ in its second argument we can
further refine this to
\[
\widetilde{\rho_L}\left(wn\right) =\sum_{y\in
  W_L}a_y\mathsf{Tr}\left(\gamma\left(wn,y\right)\right)
\] 
where the numbers $a_y$ are such that $e_L^L=\sum_{y\in W_L}a_yy$.

For fixed $n\in N_L$ we define a right action $\cdot_n$ of $W_L$ on $W_L$
by
\[
y\cdot_n z=nz^{-1}n^{-1}yz
\]
for $y,z\in W_L$. Then the stabilizer $Z_n\left(y\right)$ of $y\in W_L$
under this action is $C_{W_L}\left(n^{-1}y\right)$. We denote the orbit of
$y$ by $\mathcal{O}_n\left(y\right)=\left\{nz^{-1}n^{-1}yz\mid z\in
  W_L/C_{W_L}\left(n^{-1}y\right)\right\}$. Now
\begin{align*}
  \mathsf{Tr} \left( \gamma \left(wn,y\right)\right) &=\left|\left\{z\in
      W_L\mid w^{-1}=y\cdot_n z\right\}\right| \\
  &=\begin{cases}0&\text{if $w^{-1}\not\in\mathcal{O}_n\left(y\right)$}\\
    \left|Z_n\left(y\right)\right| &\text{if
      $w^{-1}\in\mathcal{O}_n\left(y\right)$}\end{cases} \\
  &=\begin{cases}0&\text{if $y\not\in\mathcal{O}_n\left(w^{-1}\right)$}\\
    \left|C_{W_L}\left(wn\right)\right| &\text{if
      $y\in\mathcal{O}_n\left(w^{-1}\right)$,}\end{cases}
\end{align*}
where in the last equality we have used the fact that
$w^{-1}\in\mathcal{O}_n\left(y\right)$ if and only if
$y\in\mathcal{O}_n\left(w^{-1}\right)$, and if so, then the value of
$\mathsf{Tr}\left(\gamma\left(wn,y\right)\right)$ is
$\left|Z_n\left(w^{-1}\right)\right|=\left|C_{W_L}\left(n^{-1}w^{-1}\right)\right|
=\left|C_{W_L}\left(wn\right)\right|$
by the calculation above. In
conclusion, we obtain the following formula.
\begin{align*}
  \widetilde{\rho_L}\left(wn\right)
  &=\sum_{y\in\mathcal{O}_n\left(w^{-1}\right)}a_y\left|C_{W_L}\left(wn\right)\right|\\
  &=\left|C_{W_L}\left(wn\right)\right|
  \sum_{y\in\mathcal{O}_n\left(w^{-1}\right)}
  \sum_{\substack{J\subseteq L\\\mathcal{D}\left(y\right)\subseteq L\setminus J}}n_{LJ}\\
  &=\left|C_{W_L}\left(wn\right)\right| \sum_{J\subseteq
    L}n_{LJ}\left|\mathcal{O}_n\left(w^{-1}\right) \cap X_J^L\right|.
\end{align*}
Here we have used the descent set $\mathcal{D}\left(y\right) =\left\{s\in
  L\mid\ell\left(sy\right)<\ell\left(y\right)\right\}$ to derive the formula
$\sum_{\mathcal{D}\left(y\right)\subseteq L\setminus J}n_{LJ}$ for $a_y$.

\subsection{The Extension
  \texorpdfstring{$\widetilde{\omega_L}$}{omegaLtilde}}
In this subsection we discuss the calculation of $\widetilde{\omega_L}$ for
$L\subseteq S$. As this calculation is almost identical to the calculation
of $\omega$, we begin with $\omega$ and discuss the minor modifications
needed to calculate $\widetilde{\omega_L}$ at the end.

For computational purposes, rather than working with the set $T$ of
reflections in $W$, it is simpler to work with the positive roots of $W$.
The positive roots are stored in \CHEVIE\ as vectors in the {\ttfamily
  roots} component of a Coxeter group record, the first half being the
positive roots and the second half being the negative of the first half.
This means that whenever a calculation involving roots results in a negative
root, we need to replace the negative root with its positive counterpart.

With this convention the generator $a_t$ of $A\left(W\right)$ is denoted by
$a_r$, where $r$ is the positive root orthogonal to hyperplane fixed by $t$.
To simplify the notation, we will denote $a_r$ simply by $r$. This also
reflects the way one implements $A\left(W\right)$ on a computer.  Namely,
the elements of $A\left(W\right)$ are represented by linear combinations of
sequences $r_1r_2\cdots r_q$ of positive roots.  We will also assume that
any element $r_1r_2\cdots r_q$ satisfies $r_1<r_2<\cdots<r_q$, explicitly
sorting the factors and inserting the appropriate sign $\pm 1$ whenever the
factors become unsorted. Here $<$ denotes a fixed total order on the
positive roots, which can be simply be taken to be the order in which the
roots appear in the {\ttfamily roots} component of the record for~$W$.

Now since \CHEVIE\ implements the element $w$ of $W$ as a permutation
$\sigma_w$ of the roots in $V$, it follows that if $t$ is the reflection
defined by the root $r$, then the conjugate $w^{-1}tw$ is the reflection
defined by $r.\sigma_w$, which we simplify to $r.w$. Therefore, the action
of $W$ on $A\left(W\right)$ is given by $\left(r_1r_2\cdots
  r_q\right).w=\left(r_1.w\right)\left(r_2.w\right)\cdots\left(r_q.w\right)$.

We use the non-broken circuit basis $\mathcal{B}$ of $A\left(W\right)$
described in \cite{rank3and4} to calculate its character $\omega$.  While
this works exactly as in \cite{rank3and4}, we briefly describe some
improvements to the algorithm that make the calculations in this paper
possible. Let $n=\left|S\right|$ be the rank of $W$ and recall that the
non-broken circuit basis of $A\left(W\right)$ consists of elements of the
form $r_1r_2\cdots r_n$ not containing certain sequences called {\em broken
  circuits} as subsequences. A broken circuit $r_{i_1}r_{i_2}\cdots r_{i_q}$
has the property that there exists a positive root $r$ with $r>r_{i_q}$ for
which $r_{i_1}r_{i_2}\cdots r_{i_q}r$ is dependent, so the defining relation
for $A\left(W\right)$ implies that
\begin{equation}\label{boundary}
  \left(-1\right)^q r_{i_1}r_{i_2}\cdots r_{i_q}
  =\sum_{k=1}^q\left(-1\right)^kr_{i_1}r_{i_2}\cdots \widehat{r_{i_k}}
  \cdots r_{i_q}r. 
\end{equation}
Therefore, any element {\em not} in $\mathcal{B}$ can be expressed as a
linear combination of lexicographically larger elements of $A\left(W\right)$
by applying (\ref{boundary}) to a broken circuit subsequence.  This
observation is the rationale for the procedure for expressing an arbitrary
element of $A\left(W\right)$ in terms of the non-broken circuit basis, but
it also leads to a significant improvement in the calculation of $\omega$.

Namely, to calculate the value of $\omega$ at an element $w\in W$, one in
principle runs through all basis elements $b\in\mathcal{B}$, expressing
$b.w$ as a linear combination of elements of $\mathcal{B}$ using
(\ref{boundary}) and storing the coefficients of the result into the rows of
a matrix $m$.  Then $m$ represents the linear transformation $w$ of $A\left(W\right)$
and $\omega\left(w\right)$ is the trace of $m$. We observe
that if at any point in the calculation of $b.w$ we arrive at a monomial
lexicographically larger than $b$, then this monomial cannot contribute to
the trace of $m$. Such calculations can therefore be terminated.
Furthermore, the matrix $m$ itself exists only in concept. In
practice we need only its diagonal entries. Therefore, we use the
following algorithm.

\bigskip\noindent{\bf COEFF} {\sl (Individual coefficient with respect to
  $\mathcal{B}$)} With respect to the non-broken circuit basis $\mathcal{B}$ of
$A\left(W\right)$ this algorithm takes as input a monomial $a=r_1r_2\cdots
r_n\in A\left(W\right)$ and a basis element $b\in\mathcal{B}$.  It returns
the coefficient of $b$ when $a$ is expressed with
respect to $\mathcal{B}$.

\begin{algorithmic}
  \If{$a>b$} \State\Return $0$
  \ElsIf{$a\in\mathcal{B}$}\If{$a=b$}\State\Return $1$\Else\State\Return $0$\EndIf
  \Else \State find a subsequence
  $r_{i_1}r_{i_2}\cdots r_{i_q}$ of $a$ which is a broken circuit \State find
  a root $r$ for which $r_{i_1}r_{i_2}\cdots r_{i_q}r$ is dependent \State
  \Return $\sum_{j=1}^m\left(-1\right)^j\textbf{COEFF}\left(
    r_1r_2\cdots\widehat{r_{i_j}}\cdots r_nr,b\right)$
  \EndIf
\end{algorithmic}

Observe that in the last line of the algorithm, we have inserted $r$ at the
end of the first argument of $\textbf{COEFF}$ for notational
convenience. Moving the factor to its proper position will introduce a sign
$\pm 1$.  Then to calculate $\omega\left(w\right)$ we simply calculate
$\sum_{b\in\mathcal{B}}\textbf{COEFF}\left(b.w,b\right)$.

Finally, to calculate the character $\widetilde{\omega_L}$ for $L\subseteq
S$ we calculate the non-broken circuit basis of the top component of
$A\left(W_L\right)$. Observe that an element $w\in N_W\left(W_L\right)$ is
implemented as a permutation $\sigma_w$ of the roots in $V$, so to apply $w$
to an element $r_1r_2\cdots r_q$ of $A\left(W_L\right)$ each $r_i$ must be
replaced with its corresponding root in $V$.  In \CHEVIE\ this can be
accomplished with the {\ttfamily rootInclusion} component of the $W_L$
record. Then the permutation $\sigma_w$ can be applied directly, followed by
replacing each root with the corresponding root in the reflection
representation of $W_L$ using the {\ttfamily rootRestriction} component of
the $W_L$ record. With this modification, we proceed exactly as in the
calculation of $\omega$ above.

\section{Proof of \autoref{ConjectureC} when \texorpdfstring{$L=S$}{L=S}}\label{ConjectureBSection}
Observe that if $L=S$ then $\widetilde{\rho_S}=\rho_S$,
$\widetilde{\omega_S} =\omega_S$, and $N_W\left(W_L\right)=W$. Observe also
that $\alpha_S\left(w\right)=1$ for all $w\in W$ since the space of fixed
points of $W$ is the zero subspace of $V$.  Therefore, to verify
\autoref{ConjectureC} we need to find a character $\varphi_w$ of
$C_W\left(w\right)$ for each $w\in\mathcal{C}_S$ such that
\begin{equation}\label{ConjectureB}
  \rho_S=\sum_{w\in\mathcal{C}_S} \mathrm{Ind}_{C_W(w)}^W \varphi_w =\epsilon
  \omega_S.
\end{equation}
In this section we exhibit these characters for each irreducible Coxeter
group $W$ of rank five or six. Once the characters $\varphi_w$ are
specified, one verifies (\ref{ConjectureB}) by routine calculations, so we
limit ourselves to displaying the characters $\mathrm{Ind}_{
  C_W\left(w\right)} ^W\varphi_w$ (denoted simply by $\varphi_w$),
$\epsilon$, $\rho_S$, and $\omega_S$ only for the group
$W=W\left(E_6\right)$.

Because each character $\varphi_w$ is one-dimensional, it suffices to list
its values on a generating set for the group $C_W\left(w\right)$. For the
group $W\left(E_6\right)$  we have constructed 
generating sets for the groups $C_W\left(w\right)$ ad hoc.  In type~$B$
generating sets for $C_W\left(w\right)$ are known, while in type~$D$
generating sets for $C_W\left(w\right)$ can be determined as described
below.  We use the notation for these generating sets from \cite{rank3and4},
which we now briefly review.

The cuspidal classes of $W\left(B_n\right)$ are indexed by partitions of
$n$. We always display partitions in non-decreasing order without
punctuation. With the labeling of the elements of $S$ given by the diagram
$\begin{xy}<.75cm,0cm>: (0,0)="1"; (1,0)="2"**\dir{=}; ?*{<};
  "2";(2,0)="3"**\dir{-}; (2.5,0)**\dir{-}; (3,0)="4"; (3.5,0);
  (4,0)="5"**\dir{-}; (5,0)="6"**\dir{-}; "1"*{\bullet}*+!U{_1};
  "2"*{\bullet}*+!U{_2}; "3"*{\bullet}*+!U{_3}; "4"*{\cdots};
  "5"*{\bullet}*+!U{_{n{-}1}}; "6"*{\bullet}*+!U{_n};
\end{xy}$ we define the following elements of $W\left(B_n\right)$, where we
denote the elements of $S$ by $1,2,\ldots,n$ rather than
$s_1,s_2,\ldots,s_n$ to improve legibility. If
$\lambda=\lambda_1\lambda_2\cdots\lambda_k$ is a partition of $n$ then for
each $1\le i\le k$ we define a negative $\lambda_i$-cycle
\begin{equation}\label{CIDef}
  c_i=\left(j+1\right)j\left(j-1\right)\cdots 212\cdots
  \left(j+\lambda_i\right)\qquad\text{where}\qquad
  j=\sum_{k=1}^{i-1}\lambda_k.
\end{equation}
Then each $c_i$ centralizes the element $w_\lambda=c_1c_2\cdots c_k$, which
we take to be the representative of the cuspidal class labeled by $\lambda$.
Whenever $\lambda_i=\lambda_{i+1}$ the element
\begin{equation}\label{XIDef}
  x_i=\prod_{k=1}^{\lambda_i} \left(j+k\right) \left(j+k-1\right) \cdots
  \left(j+k-\lambda_i+1\right) \qquad\text{where}\qquad
  j=\sum_{k=1}^i\lambda_k
\end{equation}
also centralizes $w_\lambda$.  In fact, if $m\left(j\right)=\min\left\{k\mid
  \lambda_k=j\right\}$ then $C_W\left(w_\lambda\right)$ is generated by the
elements $x_i$ for all $i$ satisfying $\lambda_i=\lambda_{i+1}$, together
with the elements $c_{m\left(j\right)}$ for all $j$ appearing as parts of
$\lambda$.  We remark that the elements defined in (\ref{CIDef}) and
(\ref{XIDef}) coincide with the elements $c_i$ and $x_i$ defined in
\cite{rank3and4}. The character $\varphi_{c_\lambda}$ of
$C_W\left(w_\lambda\right)$ is denoted simply by $\varphi_\lambda$.

We view $W\left(D_n\right)$ as a reflection subgroup of $W\left(B_n\right)$
generated by the reflections $1'=121$ and $2,3,\ldots,n$.  Then
$w_\lambda\in W\left(D_n\right)$ whenever $\lambda$ has an even number of
parts. In fact, such elements $w_\lambda$ are representatives of the
cuspidal classes of $W\left(D_n\right)$ and the centralizer
$C_{W\left(D_n\right)}\left(w_\lambda\right)$ is the intersection of
$W\left(D_n\right)$ with $C_{W\left(B_n\right)}\left(w_\lambda\right)$.  We
observe that the last factor $j+k-\lambda_i+1$ of (\ref{XIDef}) is at least
$\lambda_i+1$ and that the other factors are greater than
$j+k-\lambda_i+1$. This means that $1$ is never a factor of $x_i$ so that
$x_i\in W\left(D_n\right)$.  However, (\ref{CIDef}) shows that $1$ occurs as
a factor of $c_i$ exactly once, making the replacement of $121$ by $1'$
impossible. This shows that $c_i\not\in W\left(D_n\right)$.  Nevertheless,
generators of $C_{W\left(D_n\right)}\left(w_\lambda\right)$ can often be
found among products of an even number of the elements $c_i$.

In each of the following subsections we present the results of our
calculations for the finite irreducible Coxeter groups of rank five and six.
For each cuspidal class representative $w$ we display a
generating set of $C_W\left(w\right)$, where the generators are written as
words in the Coxeter generators. At each generator, we display the value of
the character $\varphi_w$.  If $\zeta$ is an eigenvalue of $w$ on $V$, we
denote the determinant of the representation of $C_W\left(w\right)$ on the
$\zeta$-eigenspace of $w$ in $V$ by $\det |_\zeta$. If $\varphi_w$ is a
power of $\det |_\zeta$ for some $\zeta$, then we also indicate this.  By
Springer's theory of regular elements \cite{springer}, the centralizer
$C_W\left(w\right)$ is a complex reflection group when $w$ is a regular
element. When this is the case, we identify $C_W\left(w\right)$ as such a
group.  For $n\ge 1$ we denote the $n^\text{th}$ root of unity $e^{2\pi
  i/n}$ by $\zeta_n$, the cyclic group of size $n$ by $Z_n$, and the
symmetric group on $n$ letters by $S_n$.

\subsection{\texorpdfstring{$W=W\left(E_6\right)$}{W=W(E6)}}
We begin with $W(E_6)$ and present the calculations that lead to the proof of
\autoref{ConjectureC} for this group. For the other groups of rank
five and six we present only the basic information described above.

Define the characters $\varphi_d= \varphi_{w_d}$ in the following table,
where the conjugacy classes of $W$ are labeled by their Carter diagrams $d$.
Here the elements of $S$ are labeled as in the Coxeter graph
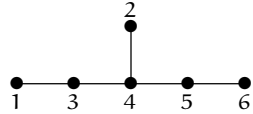
\[
\vcenter{
  \begin{xy}<.75cm,0cm>: (0,0)*{\bullet}*+!U{_1}; (1,0)*{\bullet}*+!U{_3};
    (2,0)*{\bullet}*+!U{_4}; (3,0)*{\bullet}*+!U{_5};
    (4,0)*{\bullet}*+!U{_6}; (0,0);(4,0)**\dir{-}; (2,0);(2,1)**\dir{-};
    (2,1)*{\bullet}*+!D{_2};
  \end{xy}}
\] 
and $r$ denotes the reflection defined by the highest root of $W$.

\[
\begin{array}{ccclcc}
  \toprule
  d&w_d&\text{Gen}&\varphi_d&C_W\left(w_d\right)&\text{Det}\\\otoprule
  A_2^3&12356r&24542314&\zeta_3&G_{25}&\det |_{\zeta_3}\\
  &&13&\zeta_3\\
  &&56&\zeta_3\\\midrule
  E_6(a_2)&w_{E_6}^2&w_{E_6}&1&G_5&\\
  &&234543&\zeta_3\\\midrule
  A_5A_1&13456r&w_{A_5A_1}&\zeta_3\\
  &&2345432&-1\\
  &&r&-1\\\midrule
  E_6(a_1)&34w_{E_6}&w_{E_6(a_1)}&\zeta_9&Z_9&\det |_{\zeta_9}\\\midrule
  E_6&123456&w_{E_6}&-1&Z_{12}
  &\left(\det |_{\zeta_{12}}\right)^6 \\ \bottomrule
\end{array}
\]

\medskip

Finally, the values of the characters $\varphi_d^W$ together with $\rho_S$
and $\omega_S$ are shown in the following table.

\begin{multline*} 
  \tiny \renewcommand{\arraystretch}{1.2}\arraycolsep4pt
  \begin{array}{c|ccccccccccccc}
      W&\emptyset&A_1^4&A_1^2&A_2^3&A_2&A_2^2&D_4(a_1)&A_3A_1&A_4
      &E_6(a_2)&D_4&A_5A_1&A_2A_1^2\\\midrule
      \varphi_{A_2^3}&80&16&\cdot&-10&-4&2&8&\cdot&\cdot&-2&-2&-2&\cdot\\
      \varphi_{E_6(a_2)}&720&16&\cdot&-18&-12&-6&8&\cdot&\cdot&-2&-2&-2&\cdot\\
      \varphi_{A_5A_1}&1440&32&\cdot&-36&12&-3&\cdot&\cdot&\cdot&-4&2& -1&
      \cdot\\
      \varphi_{E_6(a_1)}&5760&\cdot&\cdot&-72&\cdot&\cdot& \cdot&\cdot& \cdot&
      \cdot&\cdot&\cdot&\cdot\\
      \varphi_{E_6}&4320&96&\cdot&108&\cdot&\cdot&-16&\cdot&\cdot&12& \cdot&
      \cdot&\cdot\\\midrule
      \epsilon&1&1&1&1&1&1&1&1&1&1&1&1&1\\
      \rho_S&12320&160&\cdot&-28&-4&-7&\cdot&\cdot&\cdot&4&-2&-5&\cdot\\
      \omega_S&12320&160&\cdot&-28&-4&-7&\cdot&\cdot&\cdot&4&-2&-5&\cdot
    \end{array}\\ \\  \tiny
    \renewcommand{\arraystretch}{1.2}\arraycolsep4pt
    \begin{array}{c|cccccccccccc}
      W&E_6(a_1)&E_6&A_1&A_1^3&A_3A_1^2&A_3&A_2A_1
      &A_2^2A_1&A_5&D_5&A_4A_1&D_5(a_1)\\ \midrule
      \varphi_{A_2^3}&-1&2&\cdot&\cdot&\cdot&\cdot&\cdot&\cdot&\cdot& \cdot&
      \cdot& \cdot\\
      \varphi_{E_6(a_2)}&\cdot&2&\cdot&\cdot&\cdot&\cdot&\cdot&\cdot& \cdot&
      \cdot&\cdot&\cdot\\
      \varphi_{A_1A_5}&\cdot&\cdot&-120&-8&\cdot&\cdot&\cdot&3&1&\cdot&
      \cdot& \cdot\\
      \varphi_{E_6(a_1)}&\cdot&\cdot&\cdot&\cdot&\cdot&\cdot&\cdot&\cdot&
      \cdot& \cdot&\cdot&\cdot\\
      \varphi_{E_6}&\cdot&-4&\cdot&\cdot&\cdot&\cdot&\cdot&\cdot&\cdot&
      \cdot& \cdot&\cdot\\ \midrule
      \epsilon&1&1&-1&-1&-1&-1&-1&-1&-1&-1&-1&-1\\
      \rho_S&-1&\cdot&-120&-8&\cdot&\cdot&\cdot&3&1&\cdot&\cdot&\cdot\\
      \omega_S&-1&\cdot&120&8&\cdot&\cdot&\cdot&-3&-1&\cdot&\cdot&\cdot
    \end{array}
  \end{multline*}

%\vfill\eject

\subsection{\texorpdfstring{$W=W\left(B_5\right)$}{W=W(B5)}}
The characters defined in the following table satisfy \autoref{ConjectureC}
for $W=W\left(B_5\right)$ when $L=S$.

\[
\begin{array}{cccrcc}
  \toprule
  \lambda&\text{Gen}&\text{Word}&\varphi_\lambda
  &C_W\left(w_\lambda\right)&\text{Det}\\ \otoprule
  1^5&S&&\epsilon&W&\det |_{-1} \\  \midrule
  1^32&c_1&1&-1\\
  &c_4&43212345&-1\\
  &x_1&2&-1\\
  &x_2&3&-1\\ \midrule
  12^2&c_1&1&-1\\
  &c_2&2123&-1\\
  &x_2&4354&-1\\ \midrule
  1^23&c_1&1&-1\\
  &c_3&3212345&\zeta_6\\
  &x_1&2&-1\\ \midrule
  23&c_1&12&-1\\
  &c_2&3212345&\zeta_6\\ \midrule
  14&c_1&1&-1\\
  &c_2&212345&-1\\ \midrule
  5&c_1&12345&\zeta_{10}&Z_{10}&\det |_{\zeta_{10}}\\
  \bottomrule
\end{array}
\]

\vfill\eject

\subsection{\texorpdfstring{$W=W\left(B_6\right)$}{W=W(B6)}}
The characters defined in the following table satisfy \autoref{ConjectureC}
for $W=W\left(B_6\right)$ when $L=S$.

\[
\begin{array}{cccrcc}
  \toprule 
  \lambda&\text{Gen}&\text{Word}&\varphi_\lambda
  &C_W\left(w_\lambda\right)&\text{Det}\\\otoprule
  1^6&S&&\epsilon&W&\det |_{-1}\\ \midrule
  1^42&c_1&1&-1\\
  &c_5&5432123456&-1\\
  &x_1&2&-1\\
  &x_2&3&-1\\
  &x_3&4&-1\\\midrule
  1^22^2&c_1&1&-1\\
  &c_3&321234&-1\\
  &x_1&2&-1\\
  &x_3&5465&-1\\\midrule
  2^3&c_1&12&-1&Z_4\wr S_3&\\
  &x_1&3243&-1\\
  &x_2&5465&-1\\\midrule
  1^33&c_1&1&-1\\
  &c_4&432123456&\zeta_6\\
  &x_1&2&-1\\
  &x_2&3&-1\\\midrule
  123&c_1&1&-1\\
  &c_2&2123&-1\\
  &c_3&432123456&\zeta_6\\\midrule
  3^2&c_1&123&\zeta_6&Z_6\wr S_2&\det |_{\zeta_6}\\
  &x_1&432543654&-1\\\midrule
  1^24&c_1&1&-1\\\
  &c_2&32123456&-1\\
  &x_1&2&-1\\\midrule
  24&c_1&12&-1\\
  &c_2&32123456&-1\\\midrule
  15&c_1&1&-1\\
  &c_2&2123456&\zeta_{10}\\\midrule
  6&c_1&123456&\zeta_6&Z_{12}&\left(\det |_{\zeta_{12}}\right)^2\\\bottomrule
\end{array}
\]

\vfill\eject

\subsection{\texorpdfstring{$W=W\left(D_5\right)$}{W=W(D5)}}
The characters defined in the following table satisfy \autoref{ConjectureC}
for $W=W\left(D_5\right)$ when $L=S$.

\[
\begin{array}{cccrcc}\toprule
  \lambda&\text{Gen}&\text{Word}&\varphi_\lambda
  &C_W\left(w_\lambda\right)&\text{Det}\\\otoprule
  1^32&w_{1^32}&1'2321'3431'2345&\zeta_4\\
  &c_1&1'&-1\\
  &x_1&2&-1\\
  &x_2&3&-1\\\midrule
  23&w_{23}&1'3w_{14}&\zeta_{12}\\\midrule
  14&w_{14}&1'2345&\zeta_8&Z_8&\det |_{\zeta_8} \\ \bottomrule
\end{array}
\]

\subsection{\texorpdfstring{$W=W\left(D_6\right)$}{W=W(D6)}}
The characters defined in the following table satisfy \autoref{ConjectureC}
for $W=W\left(D_6\right)$ when $L=S$.

\[
\begin{array}{cccrcc}
  \toprule
  \lambda&\text{Gen}&\text{Word}&\varphi_\lambda
  &C_W\left(w_\lambda\right)&\text{Det}\\\otoprule
  1^6&S&&\epsilon&W&\det |_{-1} \\\midrule
  1^22^2&c_1c_3&31'234&\zeta_4\\
  &x_1&2&-1\\
  &x_3&5465&-1\\\midrule
  1^33&c_1c_4&43w_{15}&\zeta_3\\
  &x_1&2&-1\\
  &x_2&232&-1\\\midrule
  3^2&w_{3^2}&1'343w_{15}&\zeta_3&G\left(6,2,2\right)&\det |_{\zeta_6}\\
  &c_1^2&1'232&\zeta_3\\
  &x_1&432543654&-1\\\midrule
  24&w_{24}&1'3w_{15}&\zeta_8^3\\
  &c_2^2&\left(3w_{15}\right)^2&\zeta_4\\\midrule
  15&w_{15}&1'23456&\zeta_5&Z_{10}&\left(\det |_{\zeta_{10}}\right)^2\\
  \bottomrule
\end{array}
\]

\vfill\eject

\section{Proof of \autoref{ConjectureC} when \texorpdfstring{$L$}{L} is a proper subset of \texorpdfstring{$S$}{S}}
\label{ConjectureCSection}

Recall that the normalizer in $W$ of $W_L$ factors as the semidirect product
of $W_L$ and a normalizer complement \cite{howlett:normalizers}. When the
semidirect product is a direct product, $W_L$ is called {\em bulky.} It is
shown in \cite{douglasspfeifferroehrle:inductive} that \autoref{ConjectureC}
holds if either $W_L$ is bulky or the rank of $W_L$ is two or less. Also, it
is shown in \cite{douglasspfeifferroehrle:coxeter} that
\autoref{ConjectureC} holds if $W_L$ is a direct product of Coxeter groups
of type~$A$. Thus, to prove \autoref{ConjectureA} it suffices to prove
\autoref{ConjectureC} for all pairs $W,W_L$ where the rank of $W$ is five or
six and $L$ is a proper subset of $S$ for which the following
hold.
\begin{enumerate}
\item\label{NotBulky} $W_L$ is not bulky in $W$, 
\item $W_L$ has rank at least three, and
\item\label{NotA} $W_L$ is not a direct product of Coxeter groups of
  type~$A$.
\end{enumerate}
After consulting the table of bulky parabolic subgroups in \cite{rank3and4},
it remains to consider the pairs shown in~\autoref{eq:table}.

\begin{table}[h!t]
  \begin{equation*}
  \begin{array}{l|l}
    W&W_L\\\midrule
    B_5&A_2B_2\\
    B_6&A_2B_2,\; A_2B_3,\; A_3B_2,\; A_1^2B_2\\
    D_5&D_4\\
    D_6&D_4,\; D_5\\
    E_6&D_4\\
  \end{array}  
  \end{equation*}
  \caption{List of pairs $W, W_L$ to be considered for \autoref{ConjectureC}}
  \label{eq:table}
\end{table}

We consider each such pair $W,W_L$ separately in the following subsections.
For each pair we indicate representatives of the cuspidal conjugacy classes
of $W_L$, generators of the centralizers of these representatives, and
linear characters of the centralizers that satisfy the conclusion of
\autoref{ConjectureC}.  Additionally, we also give the values of
$\widetilde{\rho_L}$, $\widetilde{\omega_L}$, $\alpha_L$, and $\epsilon$ for
the pair $W(B_5), W(A_2B_2)$ and the pair $W(E_6), W(D_4)$.

In the following sections we use the symbol $w_n$ to denote a representative
of the $n^{\text{th}}$ conjugacy class of a group in the list of conjugacy
classes returned by the command {\ttfamily ConjugacyClasses} in \GAP.  We
denote the longest element of $W$ by $w_0$ and the longest element in $W_L$
by $w_L$.  As in \autoref{ConjectureBSection} the symbols $1,2,\ldots,n$
denote the elements of $S$.

\vfill\eject

\subsection{\texorpdfstring{$W=W\left(B_5\right)$}{W=W(B5)}}
As an illustration, we provide somewhat more detail in the case where
$W=W\left(B_5\right)$ and $L=\left\{1,2,4,5\right\}$.  The cuspidal
conjugacy classes in $W_L$ are represented by $w_{13}$ and $w_{15}$. The
centralizer of $w_{15}$ is $\left\langle w_{15}\right \rangle\times
\left\langle w_0\right \rangle =Z_{12}\times Z_2$. The centralizer of
$w_{13}$ is generated by $C_W\left(w_{15}\right)$ and $1$.  We define the
characters $\varphi_{13}$ and $\varphi_{15}$ by supplying their values at
these generators shown in the following table.
\[
\begin{array}{cll}\toprule
  L&\text{Type}&\text{Characters}\\\otoprule
  \left\{1,2,4,5\right\}&A_2B_2
  &\varphi_{13}:\left(w_{15},w_0,1\right)\mapsto\left(\zeta_3,1,-1\right)\\
  &&\varphi_{15}:\left(w_{15},w_0\right)\mapsto\left(\zeta_6,1\right)
  \\\bottomrule
\end{array}
\]

Then $\varphi_{13}^{N_W\left(W_L\right)} +\varphi_{15}^{N_W\left(W_L\right)}
=\widetilde{\rho_L}$.  This character is shown in \autoref{1245InB5}
together with $\widetilde{\omega_L}$, $\alpha_L$, and $\epsilon$. The
conjugacy classes of $N_W(W_L)$ are listed in the order determined by \GAP\
where $W$ is constructed using the command {\ttfamily
  W:=CoxeterGroup("B",5)} and $N_W(W_L)$ is constructed using the command
\begin{quotation} 
  {\ttfamily Normalizer(W,ReflectionSubgroup(W,[1,2,4,5]))}.
\end{quotation}
The classes are labeled by the orders of their elements and an additional
letter to distinguish them from one another.  We see that
$\widetilde{\rho_L}=\alpha_L\epsilon\widetilde{\omega_L}$ so that
\autoref{ConjectureC} holds for the pair $W,W_L$.  This completes the proof
of \autoref{ConjectureA} for $W$.

\begin{table}[h!b]
  \begin{multline*}\arraycolsep4pt 
    \begin{array}{c|rrrrrrrrrrrrrrrr}
      N_W\left(W_L\right)&1a&2a&3a&2b&2c&6a&2d&2e&6b&2f&2g&6c&2h&2i&
      6d\\\midrule 
      \widetilde{\rho_L}&6&\cdot&-3&6&\cdot&-3&-2&\cdot&1&-2&\cdot&1&-2&
      \cdot&1\\ 
      \widetilde{\omega_L}&6&\cdot&-3&6&\cdot&-3&2&\cdot&-1&2&\cdot&-1&2&
      \cdot&-1\\  
      \alpha_L&1&1&1&-1&-1&-1&1&1&1&-1&-1&-1&1&1&1\\
      \epsilon&1&-1&1&-1&1&-1&-1&1&-1&1&-1&1&-1&1&-1
    \end{array}\\ \\ \arraycolsep4pt
    \begin{array}{c|rrrrrrrrrrrrrrrr}
      N_W\left(W_L\right)&2j&2k&6e&4a&4b&12a&4c&4d&12b&2l&2m&6f&2n&2o&
      6g\\\midrule 
      \widetilde{\rho_L}&-2&\cdot&1&-2&\cdot&1&-2&\cdot&1&6&\cdot&-3&6&
      \cdot&-3\\ 
      \widetilde{\omega_L}&2&\cdot&-1&-2&\cdot&1&-2&\cdot&1&6&\cdot&-3&6&
      \cdot&-3\\ 
      \alpha_L&-1&-1&-1&1&1&1&-1&-1&-1&1&1&1&-1&-1&-1\\
      \epsilon&1&-1&1&1&-1&1&-1&1&-1&1&-1&1&-1&1&-1
    \end{array}
  \end{multline*}
  \caption{Characters of $N_W\left(W_L\right)$ where $W=W\left(B_5\right)$
    and $L=\left\{1,2,4,5\right\}$}
  \label{1245InB5}
\end{table}

\vfill\eject

\subsection{\texorpdfstring{$W=W\left(B_6\right)$}{W=W(B6)}}
The characters defined in the following table satisfy \autoref{ConjectureC}
for $W=W\left(B_6\right)$ and $W_L$ as in \autoref{eq:table}. For notational
convenience, set $M=\left\{1,2,3,4,5\right\}$ and $r$ is the reflection 
corresponding to the highest long root of $W$.

\[
\begin{array}{cll}\toprule
  L&\text{Type}&\text{Characters}\\\otoprule
  \left\{1,2,4,5\right\}&A_2B_2
  &\varphi_{13}:\left(w_{15},w_0,w_M,1\right)\mapsto\left(\zeta_3,1,1,-1
  \right)\\ 
  &&\varphi_{15}:\left(w_{15},w_0,w_M\right)\mapsto\left(\zeta_6,1,1
  \right)\\\midrule  
  \left\{1,2,3,5,6\right\}&A_3B_2 
  &\varphi_{12}:\left(1,2,3,56,w_0\right) \mapsto\left(-1,-1,-1,\zeta_3,-1
  \right)\\ 
  &&\varphi_{24}:\left(w_{24},1,w_0\right)\mapsto\left(\zeta_3^2,-1,-1
  \right)\\ 
  &&\varphi_{30}:\left(123,56,w_0\right)\mapsto\left(-\zeta_3,\zeta_3^2,-1
  \right)\\\midrule 
  \left\{1,2,4,5,6\right\}&A_2B_3
  &\varphi_{23}:\left(1,2,456,w_0\right)\mapsto\left(-1,-1,\zeta_4,-1
  \right)\\ 
  &&\varphi_{25}:\left(12,456,w_0\right)\mapsto\left(-1,\zeta_4,-1
  \right)\\\midrule 
  \left\{1,2,4,6\right\}&A_1^2B_2
  &\varphi_{12}:\left(1,2,4,6,5465,r\right)\mapsto\left(-1,-1,-1,-1,1,1
  \right)\\ 
  &&\varphi_{20}:\left(12,4,6,5465,r\right)\mapsto\left(-1,-1,-1,1,1\right)
  \\\bottomrule\end{array}\]

\bigskip

\subsection{\texorpdfstring{$W=W\left(D_5\right)$}{W=W(D5)}} 
The characters defined in the following table satisfy \autoref{ConjectureC}
for $W=W\left(D_5\right)$ and $W_L$ as in \autoref{eq:table}.

\[
\begin{array}{cll}\toprule
  L&\text{Type}&\text{Characters}\\\otoprule
  \left\{1',2,3,4\right\}&D_4
  &\varphi_{3}=\epsilon\\
  &&\varphi_{9}:\left(x_1,1'w_0\right)\mapsto\left(-1,\zeta_4\right)\\
  &&\varphi_{11}:\left(w_{11},w_0\right)\mapsto\left(\zeta_3,1\right)
  \\\bottomrule
\end{array}
\]

\bigskip

\subsection{\texorpdfstring{$W=W\left(D_6\right)$}{W=W(D6)}} 
The characters defined in the following table satisfy \autoref{ConjectureC}
for $W=W\left(D_6\right)$ and $W_L$ as in \autoref{eq:table}. For notational
convenience, set $M=\left\{1',2,3,4,5\right\}$ and $x_1=3243$.

\[
\begin{array}{cll}\toprule
  L&\text{Type}&\text{Characters}\\\otoprule
  \left\{1',2,3,4\right\}&D_4
  &\varphi_{3}:\left(1',2,3,4,6,w_M\right)\mapsto\left(-1,-1,-1,-1,1,1
  \right)\\ 
  &&\varphi_{9}:\left(6,x_1,2w_M\right)\mapsto\left(1,-1,\zeta_4\right)\\
  &&\varphi_{11}:\left(w_{11},6,w_M\right)\mapsto\left(\zeta_3,1,1
  \right)\\\midrule  \left\{1',2,3,4,5\right\}&D_5
  &\varphi_{7}:\left(w_{7},w_0,1',2,3\right)\mapsto\left(\zeta_4,-1,-1,-1,-1
  \right)\\
  &&\varphi_{15}:\left(w_{15},w_0\right)\mapsto\left(\zeta_{12},-1\right)\\
  &&\varphi_{17}:\left(w_{17},w_0\right)\mapsto\left(\zeta_8,-1\right)
  \\\bottomrule
\end{array}
\]

\vfill\eject 

\subsection{\texorpdfstring{$W=W\left(E_6\right)$}{W=W(E6)}}
Let $L=\left\{2,3,4,5\right\}$. The cuspidal conjugacy classes in $W_L$ are
represented by $w_{3}$, $w_{9}$, and $w_{11}$. The class containing $w_{11}$
is also labeled by the partition $13$ in the notation used in type $D_4$. It
is convenient to take $y_{13}=2354$ as a representative of this class
instead of $w_{11}$. The centralizer of $y_{13}$ is generated by
$y_{13},w_M,w_N$ where $M=\left\{2,3,4,5,6\right\}$ and
$N=\left\{1,2,3,4,5\right\}$. Then $W_M$ and $W_N$ both are of type
$D_5$. Notice that conjugation by $w_M$ exchanges $2$ with $3$ while
conjugation by $w_N$ exchanges $2$ with $5$.
Set $x_1 = 4354$.  Define the following
characters.

\[
\begin{array}{cll}\toprule
  L&\text{Type}&\text{Characters}\\\otoprule
  \left\{2,3,4,5\right\}&D_4
  &\varphi_{3}:\left(3,4,w_M,w_0\right)\mapsto\left(-1,-1,1,1\right)\\
  &&\varphi_{9}:\left(x_1,2w_M,243w_0\right)\mapsto\left(-1,-\zeta_4,\zeta_4
  \right)\\ 
  &&\varphi_{11}:\left(y_{13},w_M,w_N\right)\mapsto\left(\zeta_3,1,1\right)
  \\\bottomrule
\end{array}
\]

In this case, $N_W\left(W_L\right)\cong W\left(F_4\right)$.  Then using the
notation from \cite[Table~C.3]{geckpfeiffer:characters} for the irreducible
characters of $W\left(F_4\right)$ (which is identical to the notation used
in \CHEVIE), we have
\begin{align*}
  \varphi_3^{N_W\left(W_L\right)}
  &=\chi''_{\left(1,12\right)}\\
  \varphi_9^{N_W\left(W_L\right)}
  &=\chi'_{\left(6,6\right)}+\chi''_{\left(6,6\right)}\\
  \varphi_{11}^{N_W\left(W_L\right)} &=\chi''_{\left(2,4\right)}
  +\chi_{\left(9,2\right)} +\chi''_{\left(9,6\right)}
  +\chi_{\left(12,4\right)}\\
  \widetilde{\rho_L} &=\chi_{\left(1,12\right)}'' +\chi_{\left(2,4\right)}''
  +\chi_{\left(9,2\right)} +\chi_{\left(9,6\right)}''
  +\chi_{\left(6,6\right)}' +\chi_{\left(6,6\right)}''
  +\chi_{\left(12,4\right)}\\
  \widetilde{\omega_L} &=\chi_{\left(1,12\right)}'
  +\chi_{\left(2,16\right)}' +\chi_{\left(9,10\right)}
  +\chi_{\left(9,6\right)}' +\chi_{\left(6,6\right)}'
  +\chi_{\left(6,6\right)}'' +\chi_{\left(12,4\right)}.
\end{align*}
Now since $\alpha_L\epsilon=\chi_{\left(1,24\right)}$ is the sign character
of $W\left(F_4\right)$, the calculations above together with
\cite[Table~C.3]{geckpfeiffer:characters} show that the characters
$\varphi_3,\varphi_9,\varphi_{11}$ satisfy \autoref{ConjectureC} for
the pair $W,W_L$.

\bigskip 

{\bf Acknowledgments}: The authors would like to acknowledge support from
the DFG-priority program SPP1489 ``Algorithmic and Experimental Methods in
Algebra, Geometry, and Number Theory''.

\bibliographystyle{plainnat}
%\bibliography{rank567}

\end{document}